\documentclass{jhrs}
\usepackage[latin1]{inputenc}
\usepackage[all]{xy}
\usepackage{amsmath}
\usepackage{amsthm}
\usepackage{bbold}

\newtheorem*{maintheorem}{Theorem A}
\newtheorem{theorem}{Theorem}[section]
\newtheorem{definition}[theorem]{Definition}
\newtheorem{lemma}[theorem]{Lemma}
\newtheorem{example}[theorem]{Example}
\newtheorem{prop}[theorem]{Proposition}

\newcommand\labcxdxe{
  \begin{xy}
    <0.175cm,0cm>:
    0;(0.5,1)**@{-},
    0;(-0.5,1)**@{-},
    (0.5,-1);(0,0)**@{-},
    (0.5,-1);(1,0)**@{-},
    (-0.5,1);(-1.5,2)**@{-},
    (-0.5,1);(-0.83,2)**@{-},
    (-0.5,1);(-0.16,2)**@{-},
    (-0.5,1);(0.5,2)**@{-},
  \end{xy}
}
\newcommand\lcxabdxe{
  \begin{xy}
    <0.175cm,0cm>:
    0;(-1,1)**@{-},
    0;(-.33,1)**@{-},
    0;(.33,1)**@{-},
    0;(1,1)**@{-},
    (.33,1);(.83,2)**@{-},
    (.33,1);(-.16,2)**@{-},
    (.5,-1);0**@{-},
    (.5,-1);(1,0)**@{-},
  \end{xy}
}
\newcommand\lcxdxabe{
  \begin{xy}
    <0.175cm,0cm>:
    (-0.33,-1);0**@{-},
    (-0.33,-1);(.66,0)**@{-},
    (-0.33,-1);(-0.66,0)**@{-},
    (-0.33,-1);(-1.33,0)**@{-},
    0;(-.5,1)**@{-},
    0;(.5,1)**@{-},
    (-.5,1);(-1,2)**@{-},
    (-.5,1);(0,2)**@{-},
  \end{xy}
}
\newcommand\lbxacdxe{
  \begin{xy}
    <0.175cm,0cm>:
    0;(-1,1)**@{-},
    0;(-.33,1)**@{-},
    0;(.33,1)**@{-},
    0;(1,1)**@{-},
    (-.33,1);(-.83,2)**@{-},
    (-.33,1);(.16,2)**@{-},
    (.5,-1);0**@{-},
    (.5,-1);(1,0)**@{-},    
  \end{xy}
}
\newcommand\laxbcdxe{
  \begin{xy}
    <0.175cm,0cm>:
    0;(-1,1)**@{-},
    0;(-.33,1)**@{-},
    0;(.33,1)**@{-},
    0;(1,1)**@{-},
    (-1,1);(-1.5,2)**@{-},
    (-1,1);(-.5,2)**@{-},
    (.5,-1);0**@{-},
    (.5,-1);(1,0)**@{-},    
  \end{xy}
}
\newcommand\lbxcxade{
  \begin{xy}
    <0.175cm,0cm>:
    (0.33,-1);0**@{-},
    (0.33,-1);(-.66,0)**@{-},
    (0.33,-1);(0.66,0)**@{-},
    (0.33,-1);(1.33,0)**@{-},
    0;(-.5,1)**@{-},
    0;(.5,1)**@{-},
    (-.5,1);(-1,2)**@{-},
    (-.5,1);(0,2)**@{-},
  \end{xy}
}
\newcommand\laxbxcde{
  \begin{xy}
    <0.175cm,0cm>:
    (1,-1);0**@{-},
    (1,-1);(.66,0)**@{-},
    (1,-1);(1.33,0)**@{-},
    (1,-1);(2,0)**@{-},
    0;(-.5,1)**@{-},
    0;(.5,1)**@{-},
    (-.5,1);(-1,2)**@{-},
    (-.5,1);(0,2)**@{-},
  \end{xy}
}
\newcommand\rcxbxade{
  \begin{xy}
    <0.175cm,0cm>:
    (0.33,-1);(-0.66,0)**@{-},
    (0.33,-1);(0,0)**@{-},
    (0.33,-1);(0.66,0)**@{-},
    (0.33,-1);(1.33,0)**@{-},
    0;(-0.33,1)**@{-},
    0;(0.33,1)**@{-},
    (0.33,1);(0,2)**@{-},
    (0.33,1);(0.66,2)**@{-},
  \end{xy}
}
\newcommand\rcxbdexa{
  \begin{xy}
    <0.175cm,0cm>:
    0;(-1,1)**@{-},
    0;(-.33,1)**@{-},
    0;(.33,1)**@{-},
    0;(1,1)**@{-},
    (-.33,1);(-.83,2)**@{-},
    (-.33,1);(.16,2)**@{-},
    (-.5,-1);(-1,0)**@{-},
    (-.5,-1);0**@{-},
  \end{xy}
}
\newcommand\rcdexbxa{
  \begin{xy}
    <0.175cm,0cm>:
    (-.5,-1);0**@{-},
    (-.5,-1);(-1,0)**@{-},
    0;(-.5,1)**@{-},
    0;(.5,1)**@{-},
    (.5,1);(-.5,2)**@{-},
    (.5,1);(1.5,2)**@{-},
    (.5,1);(.16,2)**@{-},
    (.5,1);(.83,2)**@{-},
  \end{xy}
}
\newcommand\rdxbcexa{
  \begin{xy}
    <0.175cm,0cm>:
    0;(-1,1)**@{-},
    0;(-.33,1)**@{-},
    0;(.33,1)**@{-},
    0;(1,1)**@{-},
    (.33,1);(.83,2)**@{-},
    (.33,1);(-.16,2)**@{-},
    (-.5,-1);(-1,0)**@{-},
    (-.5,-1);0**@{-},
  \end{xy}
}
\newcommand\rdxcxabe{
  \begin{xy}
    <0.175cm,0cm>:
    (-0.33,-1);(0.66,0)**@{-},
    (-0.33,-1);(0,0)**@{-},
    (-0.33,-1);(-0.66,0)**@{-},
    (-0.33,-1);(-1.33,0)**@{-},
    0;(-0.33,1)**@{-},
    0;(0.33,1)**@{-},
    (0.33,1);(0,2)**@{-},
    (0.33,1);(0.66,2)**@{-},
  \end{xy}
}
\newcommand\rexbcdxa{
  \begin{xy}
    <0.175cm,0cm>:
    0;(-1,1)**@{-},
    0;(-.33,1)**@{-},
    0;(.33,1)**@{-},
    0;(1,1)**@{-},
    (1,1);(1.5,2)**@{-},
    (1,1);(.5,2)**@{-},
    (-.5,-1);(-1,0)**@{-},
    (-.5,-1);0**@{-},
  \end{xy}
}
\newcommand\rexdxabc{
  \begin{xy}
    <0.175cm,0cm>:
    (-1,-1);0**@{-},
    (-1,-1);(-.66,0)**@{-},
    (-1,-1);(-1.33,0)**@{-},
    (-1,-1);(-2,0)**@{-},
    0;(-0.33,1)**@{-},
    0;(0.33,1)**@{-},
    (0.33,1);(0,2)**@{-},
    (0.33,1);(0.66,2)**@{-},
  \end{xy}
}

\def\corolla{
  \xybox{
    (0,0);(1,1)**@{-},
    (0,0);(.6,1)**@{-},
    (0,0);(.2,1)**@{-},
    (0,0);(-.2,1)**@{-},
    (0,0);(-.6,1)**@{-},
    (0,0);(-1,1)**@{-},    
  }
}

\def\corbox{
  \xybox{
    (0,0);(1,0)**@{-},
    (0,0);(0,1)**@{-},
    (1,0);(1,1)**@{-},
    (0,1);(1,1)**@{-},
  }
}
\def\lass{
  \xybox{
    (1,0);(0,1)**@{-},
    (1,0);(1.25,.25)**@{-},
    (.75,.25);(1,.5)**@{-},
    (.5,.5);(.75,.75)**@{-},
    (.25,.75);(.5,1)**@{-},
  }
}
\def\rass{
  \xybox{
    (-1,0);(0,1)**@{-},
    (-1,0);(-1.25,.25)**@{-},
    (-.75,.25);(-1,.5)**@{-},
    (-.5,.5);(-.75,.75)**@{-},
    (-.25,.75);(-.5,1)**@{-},
  }
}

\begin{document}

\title{A partial $A_\infty$-structure on the cohomology of $C_n\times
  C_m$}

\author{Mikael Vejdemo-Johansson} 
\email{mik@math.uni-jena.de}
\address{Lehrstuhl Algebra \\
  Mathematisches Institut \\
  Fakultät für Mathematik und Informatik \\
  Friedrich-Schiller-Universität Jena \\
  07737 Jena \\
  Germany} 

\thanks{The author was partially supported by DFG Sachbehilfe grant GR
  1585/4-1.}


\received{July 11, 2007}
\revised{February 29, 2008}
\published{March 1, 2008}
\submitted{James Stasheff}

\volumeyear{2008}
\volumenumber{3}
\issuenumber{1}

\startpage{1}

\begin{abstract}
  Suppose $k$ is a finite field, and $n,m\geq 4$ multiples of the field
  characteristic. Then the $A_\infty$-structure of the group
  cohomology algebras $H^*(C_n,k)$ and $H^*(C_m,k)$ are well known. We
  give results characterizing an $A_\infty$-structure on
  $H^*(C_n\times C_m,k)$ including limits on non-vanishing low-arity
  operations and an infinite family of non-vanishing higher
  operations.
\end{abstract}

\maketitle

\section{Introduction}
\label{sec:introduction}

Let $k$ be a field of positive characteristic $p$, and $C_n$ and $C_m$
cyclic groups of order at least $4$ such that $p|n$ and $p|m$. The
ring structure of the graded commutative cohomology rings
$H^*(C_n)=H^*(C_m)=\Lambda(x)\otimes k[y]$ and $H^*(C_n\times
C_m)=H^*(C_n)\otimes H^*(C_m)=\Lambda(x,z)\otimes k[y,w]$ are well
known. However, to date there is only one family of examples in the
literature where an $A_\infty$-structure on a group cohomology algebra
has been completely described: The complete calculation of an
$A_\infty$-algebra structure on $H^*(C_n)$ was performed by Dag Madsen
in \cite{madsen_phd}.

The example of the cohomology ring of cyclic groups occurs in other
nearby fields -- Ainhoa Berciano studies it in the context of tensor
factors of $H_*(K(\mathbb Z,n); \mathbb Z_p)$, where the duals of
these cohomology rings occur with grade shifts in the generators of
the group ring.
\cite{berciano-2007}

By use of the diagonal on the associahedron, described by Samson
Saneblidze and Ron Umble in \cite{saneblidze-umble-2004}, the
$A_\infty$-structures of the cyclic group cohomologies can be extended
to $A_\infty$-structures on any finite abelian group. The exact form
these take, though, depends heavily on the actual combinatorial
details of the Saneblidze-Umble diagonal and its iterates.

The applications to group cohomology follow from a slightly more
general result, which forms the main result of this paper. 

\begin{maintheorem}
  Let $n\geq m>3$ and let $A$ and $B$ be $A_\infty$-algebras with
  $m_2\neq 0$, $m_n\neq 0$ and $m_r=0$ for all other values of $1\leq
  r<n+m$ in $A$ and $m_2\neq 0$, $m_m\neq 0$ and $m_r=0$ for all other
  values of $1\leq r<n+m$ in $B$.

  Then the only possible arities of non-trivial operations of
  $A\otimes B$ of arity less than $n+m$ are $2$, $n$, $m$
  and $n+m-2$. The operations of arity $2,n,m$ are nontrivial
  regardless of further structure on $A$ and $B$.

  Suppose finally that $n,m\geq 4$ are both divisible by $p$. Then
  $A=H^*(C_n)$ and $B=H^*(C_m)$ are non-trivial, non-formal
  $A_\infty$-algebras and all operations on $H^*(C_n\times C_m)$ of
  arities $k(n-2)+k(m-2)+2)$, $k(n-2)+(k-1)(m-2)+2$ and
  $(k-1)(n-2)+k(m-2)+2$, for $k\geq 0$ are non-trivial.
\end{maintheorem}

The paper is organized as follows: Section \ref{sec:Aooalgebra}
recalls the notion of an $A_\infty$-algebra, and gives the information
we need about the cohomology of cyclic groups. Section
\ref{sec:saneblidze-umble} recalls the construction of the
Saneblidze-Umble diagonal. Section \ref{sec:combinatorics} contains
combinatorial observations on the diagonal, and section
\ref{sec:consequences} collates the result to statements on the
cohomology ring $H^*(C_n\times C_m)$.

\section{$A_\infty$-algebras}
\label{sec:Aooalgebra}

A graded $k$-vector space $A$ is an $A_\infty$-algebra if one of the
following equivalent conditions hold

\begin{enumerate}
\item There is a family of maps $\mu_i\colon A^{\otimes i}\to A$,
  called higher multiplications fulfilling the Stasheff identities 
  \[
  \operatorname{St}_n\colon \sum_i\sum_j \mu_i\circ_j\mu_{n-i} = 0\quad.
  \]
\item There is a family of chain maps from the cellular chain complex
  of the associahedra to appropriate higher endomorphisms of $A$
  \[
  \mu_n\colon C_*(K_n)\to\operatorname{Hom}(A^{\otimes n}, A)\quad.
  \]
\item $A$ is a representation of the free dg-operad
  resolution $\mathcal Ass_\infty$ of the associative operad. 
\end{enumerate}

The structure was introduced by Jim Stasheff in
\cite{stasheff_ha_hspaces63}, and a deeper discussion suitable for the
representation theoretic point of view can be found in Bernhard
Keller's papers
\cite{keller_intro01} and in \cite{keller_Aoo_repr00} as well as in
the papers 
\cite{lu_wu_palmieri_zhang_Aoo_ring04} and
\cite{lu_wu_palmieri_zhang_Aoo_Ext06} by Lu, Palmieri, Wu and Zhang.

By a theorem by Tornike Kadeishvili \cite{kadeishvili80} and several
others, we can construct an $A_\infty$-algebra structure on the
homology $HA$ of a dg-algebra $A$ together with a quasiisomorphism of
$A_\infty$-algebras $HA\to A$.

In group cohomology, we consider $H^*(G)=\operatorname{Ext}_{kG}^*(k,k)$,
which we calculate as the homology of the endomorphism dg-algebra
$\operatorname{End}(Pk,Pk)$ of a projective resolution $Pk$ of the trivial
$kG$-module $k$. This endomorphism dg-algebra thus induces an
$A_\infty$-structure on $H^*(G)$.

Suppose $G$ is a $p$-group or an abelian group. Then the
$A_\infty$-structure on $H^*(G)$ is enough to reconstruct $kG$ up to
isomorphism, by theorems by Keller \cite{keller_intro01} and Lu,
Palmieri, Wu and Zhang \cite{lu_wu_palmieri_zhang_Aoo_Ext06}.

Johannes Huebschmann has with great success used $A_\infty$-algebra
and module structures to compute free resolutions and group cohomology
rings. These computations still give more explicit descriptions of
specific cohomology ring structures than other methods available for
computing cohomology rings.
\cite{huebschmann-coho-89,huebschmann-modp-89,huebschmann-perturbation-89,huebschmann-91}
However, in this paper we consider explicit computation of
$A_\infty$-algebra structure on the resulting group cohomology
rings. While the papers by Huebschmann certainly adress the
multiplicative structure of group cohomology rings, they do not adress
the computation of higher multiplicative structures. 
Though the existence of these structures has been known for a long
time, the actual structures are largely
uncomputed. \cite{huebschmann-personal-08} The one exception is a
structure on $H^*(C_n)=k[x,y]/(x^2)$ for appropriate cyclic groups
$C_n$ which was computed by Dag Madsen. This structure has the cup
product as $\mu_2$ and
$\mu_n(xy^{i_1},\dots,xy^{i_n})=y^{i_1+\dots+i_n+1}$. See the appendix
of \cite{madsen_phd} for details of this calculation.

\section{The Saneblidze-Umble diagonal}
\label{sec:saneblidze-umble}

Let us review the enumeration of the Saneblidze-Umble diagonal on the
cellular chains of the associahedron. This exposition follows that of
Ainhoa Berciano in \cite{berciano-2007}. For details, please refer to
\cite{saneblidze-umble-2004} or \cite{berciano-2007}.

\begin{definition}
  A \emph{step matrix} is a matrix whose non-zero entries
  \begin{itemize}
  \item Include each integer in $[n]=\{1,2,\dots,n\}$ precisely once.
  \item Occur adjacently in each row and each column.
  \item Occur strictly increasing to the right and downwards.
  \item Occur exactly once in each diagonal parallell to the main diagonal.
  \end{itemize}
\end{definition}

\begin{prop}
  Step matrices with entries from $[m]$ correspond
  bijectively to permutations in $S_m$.
\end{prop}

Next, we define right-shift and down-shift matrix transformations.

\begin{definition}
  Given a $r\times s$-matrix $G=(g_{i,j})$, we define
  \begin{itemize}
  \item for $M_j$ a non-empty subset of the non-zero entries in column
    $j$, $R_{M_j}G$ is the matrix interchanging each $g_{k,j}\in M_j$
    with $g_{k,j+1}$ if
    \begin{itemize}
    \item $\min M_j > \max\{g_{*,j+1}\}$ and
    \item $g_{t,j+1}=0$ for $g_{k,j}=\min M_j$ and $k\leq t\leq r$.
    \end{itemize}
    otherwise, define $R_{M_j}G=G$.
  \item for $N_j$ a non-empty subset of the non-zero entries in row
    $j$, $D_{N_j}G$ is the matrix interchanging each $g_{j,k}\in N_j$
    with $g_{j+1,k}$ if
    \begin{itemize}
    \item $\min N_j > \max\{g_{j+1,*}\}$ and
    \item $g_{j+1,t}=0$ for $g_{j,k}=\min N_j$ and $k\leq t\leq s$.
    \end{itemize}
    otherwise define $D_{N_j}G=G$.
  \end{itemize}
\end{definition}

\begin{definition}
  Suppose $G$ is a step matrix. Then a \emph{derived matrix}, derived
  from $G$, is a matrix of the form
  \[
  D_{N_i}D_{N_{i-1}}\dots D_{N_1}R_{M_j}R_{M_{j-1}}\dots R_{M_1} G.
  \]
\end{definition}

Note that step matrices are derived matrices via $N_i=\emptyset$,
$M_j=\emptyset$ for all $i,j$. 

\begin{definition}
  Let $\lambda_A=A_1|A_2|\dots|A_s$ and $\lambda_B=B_1|B_2|\dots|B_r$
  be partitions of $[n]$. We call the pairing $\lambda_A\otimes
  \lambda_B$ an \emph{$(s,r)$-complementary pairing} (CP) if there is
  an $r\times s$ derived matrix with columns $A_1,\dots,A_s$ and rows
  $B_r,\dots,B_1$.
\end{definition}

Complementary pairings correspond in an obvious way bijectively with
derived matrices. Partitions of $[n]$ in turn correspond to planar
rooted leveled trees with $\lambda_A=A_1|\dots|A_s$ corresponding to a
tree with root in level $s$, $n+1$ leaves, each $A_i$ describing the
corollas in level $i$ with $j\in A_i$ indicating that the branch
containing the leaf $j$ will meet the branch containing the leaf
$j+1$ in the level $i$.

Using this correspondence, we can now define a diagonal on the
permutahedron.

\begin{definition}
  Denoting the top dimensional cell of $P_n$ by $e^n$, we define
  $\Delta_P(e^0)=e^0\otimes e^0$. Inductively, having defined
  $\Delta_P$ on $C_*(P_{k+1})$ for all $0\leq k\leq n-1$, we define
  $\Delta_P$ on $C_n(P_{n+1})$ by 
  \[
  \Delta_P(e^n)=\sum u\otimes v
  \]
  where the sum is taken over all $(s,r)$-complementary pairings
  $u\otimes v$ with $s+r=n+2$, and we extend multiplicatively to all
  of $C_*(P_{n+1})$.
\end{definition}

The faces of the permutahedron are indexed by these planar rooted
leveled trees. To obtain cellular chains on the associahedron we apply
the projection from \cite{tonks97}, which on a tree level forgets
about the levels. When we do this, however, we will get degenerate
faces, characterized by having several corollas on the same level.

Let $\theta\colon C_*(P_*)\to C_*(K_*)$ to be the Tonks
projection to the associahedron. Degenerate faces will map to $0$ for
dimensional reasons. Using this, we can define the Saneblidze-Umble
diagonal $\Delta_K$.

\begin{definition}
  $\Delta_K\colon C_*(K_{n+2})\to C_*(K_{n+2})\otimes C_*(K_{n+2})$
  is defined by
  \[
  \Delta_K\theta = (\theta\otimes\theta)\Delta_P\quad.
  \]
\end{definition}

\section{Combinatorics on the diagonal}
\label{sec:combinatorics}

The matrices that arise in the definition of the Saneblidze-Umble
diagonal relay a lot of information about the tree structures in the
various terms of the diagonal. Most of the combinatorial background
here is known to Ainhoa Berciano and Ron Umble
\cite{umble_personal07}, but has not yet appeared in published
form. Hence, for completeness, we give the relevant statements and
their justification here. 

\begin{definition}
  We say that two entries $g_{i,j},g_{i+1,j}$ in column $j$ of a
  derived matrix $G$ are \emph{derived consecutive}, if all $k$ in the
  range $g_{i,j}<k<g_{i+1,j}$ occur in columns further left in the
  matrix. We say, dually, that two entries $g_{i,j},g_{i,j+1}$ in a
  row $i$ of a derived matrix $G$ are derived consecutive if all $k$
  in the range $g_{i,j}<k<g_{i+1,j}$ occur in rows lower down in the
  matrix.
\end{definition}

\begin{lemma}\label{lemma:blockcorolla}
  Each column of a derived matrix divides into derived
  consecutive blocks whose lengths index the
  orders of the corollas that will appear in that level.
\end{lemma}
\begin{proof}
  Suppose $a_1,\dots,a_m$ are derived consecutive in row or column
  $j$. Then the levels preceeding $j$ in the graph will have already
  connected all $a_i+1,\dots,a_{i+1}$, for all the elements failing to
  appear in the sequence $a_1,\dots,a_m$. Thus, in order for all $a_i$
  to meet $a_i+1$ at the level $j$, all the subtrees already
  connecting all the gaps have to meet in one single corolla. Thus,
  the derived consecutive block indexes a single corolla of arity $m+1$.
\end{proof}

\begin{lemma}
  If one factor of a term of the diagonal is constructed using only
  $m_2$, then the other factor has to be a single corolla of the
  appropriate arity.
\end{lemma}
\begin{proof}
  The proof is symmetric for the two possible locations for the
  factors, so we shall consider the case where the left factor has all
  $m_2$. This is given by the one-by-one matrix
  \[
  \begin{pmatrix}
    1 & 2 & \dots & n-1 \\ 
  \end{pmatrix}
  \]
  which has a single row which is a derived consecutive block in its
  own right, proving the claim.
\end{proof}

\begin{lemma}
  The non-degenerate terms of $\Delta_K(\theta(e^n))$ are given by
  matrices with exactly one derived consecutive block in each row and
  column.
\end{lemma}
\begin{proof}
  The proof is a direct application of lemma \ref{lemma:blockcorolla}.

  Suppose some row or column would have two disjoint derived
  consecutive blocks. In that case, there would be two or more
  corollas occuring on that level. However, this would imply that the
  face described by this matrix is degenerate, and thus vanishes.
\end{proof}

The following theorem extends and complements results obtained by
Ainhoa Berciano and Ron Umble \cite{umble_personal07} independently
of the author. Their results deal exclusively with the existences and
non-trivialities of operations of arity less than or equal to $2p-2$
in the $A_\infty$-coalgebraic case. My proof of these low-arity cases is
very similar to the arguments used by Berciano-Umble. However, the
extension of their results proving non-triviality of higher operations
of the stated arities is new.

\begin{maintheorem}\label{thm:main}
  Let $n\geq m>3$ and let $A$ and $B$ be $A_\infty$-algebras with
  $m_2\neq 0$, $m_n\neq 0$ and $m_r=0$ for all other values of $1\leq
  r<n+m$ in $A$ and $m_2\neq 0$, $m_m\neq 0$ and $m_r=0$ for all other
  values of $1\leq r<n+m$ in $B$.

  Then the only possible arities of non-trivial operations of
  $A\otimes B$ of arity less than $n+m$ are $2$, $n$, $m$
  and $n+m-2$. The operations of arity $2,n,m$ are nontrivial
  regardless of further structure on $A$ and $B$.

  Suppose finally that $n,m\geq 4$ are both divisible by $p$. Then
  $A=H^*(C_n)$ and $B=H^*(C_m)$ are non-trivial, non-formal
  $A_\infty$-algebras and all operations on $H^*(C_n\times C_m)$ of
  arities $k(n-2)+k(m-2)+2)$, $k(n-2)+(k-1)(m-2)+2$ and
  $(k-1)(n-2)+k(m-2)+2$, for $k\geq 0$ are non-trivial.
\end{maintheorem}

We shall prove the main theorem, by proving each atomic statement as a
separate lemma. This will proceed as follows: in the lemmata
\ref{lem:an-2am-2+2}, \ref{lem:anam2} and \ref{lem:anam3}, we
demonstrate the existence of specific diagonal terms of a particularly
good form. Then, in \ref{lem:nontrivial}, we demonstrate one argument
to the higher operations that vanishes on all diagonal terms not of
the form in the preceeding lemmata.

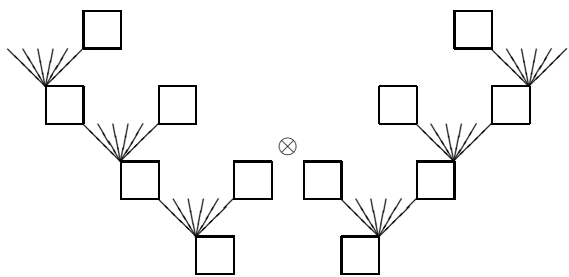
\begin{lemma}\label{lem:an-2am-2+2}
  There are diagonal terms of arity $k(n-2)+k(m-2)+2$.
\end{lemma}
\begin{proof}
  The case for $k=0$ is taken care of by the matrix 
  \[
  \begin{pmatrix}
    1
  \end{pmatrix}
  \]

  There is a derived matrix of the form 
  \begin{equation}\label{eq:derived1}
  \left(
    \begin{xy}
    <0.5cm,0cm>:
    (-2,2);
    (-2,1)**@{-};
    (-1,1)**@{-};
    (-1,0)**@{-};
    (0,0)**@{-}; 
    (0,-1)**@{-};
    (1,-1)**@{-};
    (1,-2)**@{-};
    (2,-2)**@{-};
    \end{xy}
  \right)
  \end{equation}
  where the picture is taken to depict a sparse matrix with non-zero
  entries only along the polygonal path, each horizontal line
  corresponding to $n-1$ consecutive integers and each vertical line
  corresponding to $m-1$ consecutive integers. This matrix exists
  since it can be constructed from a $k(m-2)+1\times k(n-2)+1$-matrix
  of the form
  
  \begin{equation}
  \left(
  \begin{xy}
    <0.5cm,0cm>:
    (-2,-2);
    (-2,2)**@{-};
    (2,2)**@{-};
  \end{xy}
  \right)\label{eq:step1}
  \end{equation}
  where again the polygonal path depicts the only positions in the
  matrix with non-zero entries. The sequence of moves constructing the
  matrix \eqref{eq:derived1} from the matrix \eqref{eq:step1} would
  use right shifts and down shifts that places each block in the
  zigzag where it belongs. The column in this step matrix would be a
  sequence of blocks of subsequent integers, each block of length
  $n-1$ and each block ending with an element on the form
  $k(n-2)+(k-1)(m-2)+1$. The row would start with $1$ in the first
  column, and then have a sequence of blocks of subsequent integers,
  each of length $m-1$, and each ending with an element on the form
  $k(n-2)+k(m-2)+1$.

  This matrix can be transformed into the snake like matrix given
  earlier by moving each block down or right to the expected position
  using down shifts and right shifts. Since any element that gets
  moved will move past only elements that are smaller than itself, and
  that have stopped higher up, and higher to the left, all moves
  needed are admissible.

  All in all, if we have $k$ blocks down and $k$ blocks to the right,
  the last element is $k(n-2)+k(m-2)+1$. Thus, the thus
  described operation has arity $k(n-2)+k(m-2)+2$.
\end{proof}

\begin{lemma}\label{lem:anam2}
  There are diagonal terms of arity $k(n-2)+(k-1)(m-2)+2$
\end{lemma}
\begin{proof}
  Similarily to in lemma \ref{lem:an-2am-2+2}, we can construct a
  derived matrix of the form 
  \[
  \left(
    \begin{xy}
    <0.5cm,0cm>:
    (-2,2);
    (-2,1)**@{-};
    (-1,1)**@{-};
    (-1,0)**@{-};
    (0,0)**@{-}; 
    (0,-1)**@{-};
    (1,-1)**@{-};
    (1,-2)**@{-};
    \end{xy}
  \right)
  \]
  by simply dropping the last block in the top row, and proceeding
  with everything else just as in the proof of lemma
  \ref{lem:an-2am-2+2}. The result has highest element
  $k(n-2)+(k-1)(m-2)+1$, and so the corresponding operation has arity
  $k(n-2)+(k-1)(m-2)+2$.
\end{proof}

\begin{lemma}\label{lem:anam3}
  There are diagonal terms of arity $(k-1)(n-2)+k(m-2)+2$.
\end{lemma}
\begin{proof}
  Again, similar to lemma \ref{lem:an-2am-2+2}, we can construct a
  derived matrix of the form 
  \[
  \left(
    \begin{xy}
    <0.5cm,0cm>:
    (-2,1);
    (-1,1)**@{-};
    (-1,0)**@{-};
    (0,0)**@{-}; 
    (0,-1)**@{-};
    (1,-1)**@{-};
    (1,-2)**@{-};
    (2,-2)**@{-};
    \end{xy}
  \right)
  \]
  which results from down shifts and right shifts from a matrix on the
  form 
  \[
  \left(
  \begin{xy}
    <0.5cm,0cm>:
    (-2,-2);
    (-2,2)**@{-};
    (2,2)**@{-};
  \end{xy}
  \right)
  \]
  where the first row is a sequence of blocks of subsequent integers,
  each block of length $m-1$, and each block ending with an entry on
  the form $(k-1)(n-2)+k(m-2)+1$, and the first column has a $1$ in
  the first row, and thereafter is a sequence of blocks, each of
  length $n-1$, and each ending with an entry on the form
  $k(n-2)+k(m-2)+1$. 

  This matrix has highest entry $(a-1)(n-2)+a(m-2)+1$, and so the
  corresponding operation has arity $(a-1)(n-2)+a(m-2)+2$. 
\end{proof}

\begin{lemma}\label{lem:nontrivial}
  The ``snake-like'' diagonal terms displayed above do not vanish as
  operations on $H^*(C_n\times C_m,\mathbb F_2)$.
\end{lemma}
\begin{proof}
  We shall prove the statement for the snake-like operation of arity
  $k(n-2)+k(m-2)+2$. The other two cases follow by removing runs of
  $1\otimes x$ or $x\otimes 1$ from the proposed argument, and in the
  term diagram by adding boxes to the left of the uppermost corolla on
  the left hand side or to the right of the uppermost corolla on the
  right hand side.

  First off, $H^*(C_n\times C_m,\mathbb F_2)$ has algebra generators
  $x\otimes 1$ and $1\otimes x$ of degree 1 and $y\otimes 1$ and
  $1\otimes y$ of degree 2.

  Now, we consider the argument 
  
  \begin{multline*}
    x\otimes 1, \overset{n-2\text{ times}}\cdots, x\otimes 1, x\otimes
    x, x\otimes x, 1\otimes x, \overset{m-4\text{ times}}\cdots,
    1\otimes x, x\otimes x, \\ 
    x\otimes 1,\overset{n-4\text{ times}}\cdots,x\otimes 1
    x\otimes x, \dots x\otimes x, x\otimes x,
    1\otimes x, \overset{m-2\text{ times}}\cdots, 1\otimes x\quad.
  \end{multline*}

  For a diagonal term not to vanish with this argument, it will need
  to have the form
  \[
  \begin{xy}
    <.5cm,0cm>:
    (0,-2)*{\corolla},
    (1.5,-1)*{\corbox},
    (-1.5,-1)*{\corbox},
    (.5,-3)*{\corbox},
    (-2,0)*{\corolla},
    (-.5,1)*{\corbox},
    (-3.5,1)*{\corbox},
    (-4,2)*{\corolla},
    (-2.5,3)*{\corbox},
  \end{xy}
  \otimes
  \begin{xy}
    <.5cm,0cm>:
    (0,-2)*{\corolla},
    (-1.5,-1)*{\corbox},
    (1.5,-1)*{\corbox},
    (-.5,-3)*{\corbox},
    (2,0)*{\corolla},
    (.5,1)*{\corbox},
    (3.5,1)*{\corbox},
    (4,2)*{\corolla},
    (2.5,3)*{\corbox},
  \end{xy}
  \]
  with the boxes consisting of trees built out of
  $m_2$'s, and the tree above and below each higher corolla
  containing, together, 2 less inputs than the corollas on the other
  side of the tensor product, since the running blocks of $x$'s need
  to hit the larger corollas, and the $1$'s cannot hit the larger
  corollas, lest the term vanishes.

  Thus, by considering the structure of the left hand tree, the first
  column must contain $1,2,\dots,n-2,k$, where $k$ is one more than the
  highest occuring digit in the first box. Thus, in order for the term
  not to vanish under the Tonks projection, we need $k=n-1$.

  Continuing down the tree, we get, since $k=n-1$, that after the
  column with $1,2,\dots,n-1$, we get a sequence of columns containing
  one digit each, ending with $n-1+m-2$. Then,
  $(n-2)+(m-2)+1,\dots,(n-2)+(m-2)+(n-2)$ have to occur in a single
  column, to accomodate the next corolla, and again, in order for the
  term not to vanish under the Tonks' projection, we cannot have
  anything in the box above and to the right of the corolla.

  We can continue this argument to conclude that on the left hand
  side, all the upper right boxes actually vanish.

  By symmetry, and by repeating the argument for the right hand tree
  from the bottom up, we get that all the upper left boxes vanish.

  Thus, any tree that does not vanish on the given arguments has the
  form
  
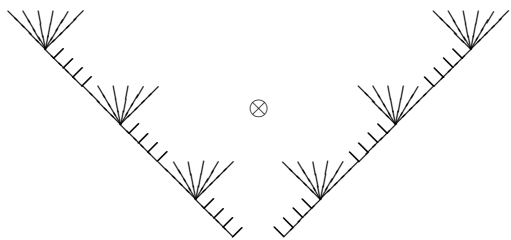
\[
  \begin{xy}
    <.5cm,0cm>:
    (.5,-3)*{\lass},
    (-0.1,-2)*{\corolla},
    (-1.5,-1)*{\lass},
    (-2.1,0)*{\corolla},
    (-3.5,1)*{\lass},
    (-4.1,2)*{\corolla},
  \end{xy}
  \otimes
  \begin{xy}
    <.5cm,0cm>:
    (-.5,-3)*{\rass},
    (0.1,-2)*{\corolla},
    (1.5,-1)*{\rass},
    (2.1,0)*{\corolla},
    (3.5,1)*{\rass},
    (4.1,2)*{\corolla},
  \end{xy}
  \]
  and is a realization of the snake-like term in lemma
  \ref{lem:an-2am-2+2}. Since any other tree pair of the same arity
  will vanish on the given arguments, this term is the only term in
  the entire diagonal sum that influences the value of the diagonal at
  this point. Hence for this particular argument we do get a
  non-vanishing value in arity $k(n-2)+k(m-2)+2$.
\end{proof}

\begin{lemma}
  If $k<n+m-1$, then either $k\in\{2,n,m,n+m-2\}$ or $m_k=0$.
\end{lemma}
\begin{proof}
  This argument was discovered independently by Ron Umble and Ainhoa
  Berciano \cite{umble_personal07}.

  $m_1$ vanishes by the way we construct the $A_\infty$-structures on
  cohomology rings.

  For the $2$-ary operation, the diagonal expression is
  $\zeta_A\otimes\zeta_B\circ\Delta_K\circ\theta(m_2)=m_2\otimes m_2$,
  which concludes the description.

  All operations of arity between $2$ and $m$ will have terms of the
  diagonal involving $m_2$ and at least one higher corolla. All higher
  corollas of arity less than $m$ vanish by the properties of the
  individual $A_\infty$-algebras.

  For arities $m$ and $n$, we have the highest order corolla
  available, and thus we will have, writing $m_2^{(n)_L}$ for a
  left-associating tower of $m_2$, and $m_2^{(n)_R}$ for a
  right-associating tower of arity $n$, the summands $m_n\otimes
  m_2^{(n)_R}$ and $m_2^{(m)_L}\otimes m_m$ non-vanishing. All other
  terms will vanish in all degrees up to the first degree in which we
  can find a non-trivial term involving both $m$-ary and $n$-ary
  operations.

  For a $k$-ary operation to contain both $m$-ary and $n$-ary
  operations, we need to fit at least one column with $n-1$ derived
  consecutive entries and one row with $m-1$ derived consecutive
  entries using only $k-1$ entries. We could share one entry between
  row and column, which gives us $m-1+n-1-1=n+m-3$ entries needed,
  which would give us a $n+m-2$-ary operation. Thus, all non-listed
  higher multiplications of arity less than $n+m-2$ will vanish.
\end{proof}

Note that if $n$ or $m$ is equal to $3$, then the snake-like terms
described above correspond to tree pairs with sequences of
higher operations composed directly without intervening
$m_2$-operations. However, by the $A_\infty$-structure on
$H^*(C_n)$, the output of a higher operation is always an even degree
element, whereas the only arguments that yield a non-vanishing higher
operations are of odd degree. Hence, such a composition will vanish.

\section{Consequences}
\label{sec:consequences}

Using Theorem A together with Madsen's results on the
$A_\infty$-structure of $H^*(C_n)$, we get a description of the
low-arity part of any non-trivial $H^*(C_n\times C_m)$.
\begin{example}
  Consider $G=C_4\times C_4$. The cohomology ring has algebra
  structure $k[x_1,x_2,y_1,y_2]/(x_1^2,x_2^2)$, and the nonzero higher
  operations involving at most $7$ arguments are given by
  \begin{align*}
    m_4(x_1a_1, x_1a_2, x_1a_3, x_1a_4) &= y_1a_1a_2a_3a_4 \\
    m_4(x_2a_1, x_2a_2, x_2a_3, x_2a_4) &= y_2a_1a_2a_3a_4 \\
    m_4(x_1x_2a_1, x_1a_2, x_1a_3, x_1a_4) &= x_2y_1a_1a_2a_3a_4 \\
    m_4(x_1x_2a_1, x_2a_2, x_2a_3, x_2a_4) &= x_1y_2a_1a_2a_3a_4 \\
    m_4(x_1a_1, x_1x_2a_2, x_1a_3, x_1a_4) &= x_2y_1a_1a_2a_3a_4 \\
    m_4(x_2a_1, x_1x_2a_2, x_2a_3, x_2a_4) &= x_1y_2a_1a_2a_3a_4 \\
    m_4(x_1a_1, x_1a_2, x_1x_2a_3, x_1a_4) &= x_2y_1a_1a_2a_3a_4 \\
    m_4(x_2a_1, x_2a_2, x_1x_2a_3, x_2a_4) &= x_1y_2a_1a_2a_3a_4 \\
    m_4(x_1a_1, x_1a_2, x_1a_3, x_1x_2a_4) &= x_2y_1a_1a_2a_3a_4 \\
    m_4(x_2a_1, x_2a_2, x_2a_3, x_1x_2a_4) &= x_1y_2a_1a_2a_3a_4
  \end{align*}
  where all the $a_i$ are monomials in $y_1,y_2$. We further get $102$
  cases of varying input configurations for $m_6$, where the deciding
  factor is that in the terms of the sum
  \[
\labcxdxe\otimes\rcxbxade+
\labcxdxe\otimes\rcxbdexa+
\labcxdxe\otimes\rcdexbxa+
\lcxabdxe\otimes\rcdexbxa+
\lcxdxabe\otimes\rcdexbxa+
\lbxacdxe\otimes\rdxbcexa+
\laxbcdxe\otimes\rdxcxabe+
\lbxcxade\otimes\rexbcdxa+
\laxbxcde\otimes\rexdxabc
  \]
  each corolla of degree 2 can see at most one odd input, and each
  corolla of degree 4 needs to see only odd input, and delivers only
  even input. Thus, for instance, the last term will have value
  $y_1y_2$, when operating on $x_2\otimes x_2\otimes x_1x_2\otimes
  x_1x_2\otimes x_1\otimes x_1$.
\end{example}

We note furthermore that by a result from Ainhoa Berciano
\cite{berciano-real-2007,berciano-2006}, the $A_\infty$-coalgebra
structure dual to the one we consider for $H^*(C_p\times C_p)$ will
have vanishing higher comultiplication in all arities except for
$q=i(p-2)+2$. We note that all such $i(p-2)+2$ occur as
$k(p-2)+k(p-2)+2$ or $k(p-2)+(k+1)(p-2)+2$ and therefore in this dual
case, any possibly non-vanishing structure map is actually
non-vanishing.

\section{Acknowledgments}
\label{sec:acknowledgments}

Many of the ideas in this paper matured shortly after Ron Umble
explained the construction of the diagonal to the author, and he is
worth deep thanks for the continuous stream of comments and
encouragement since. The paper would not be here without those.

The author is grateful for the kind and helpful suggestions put forth
by the referee and for the resulting improvement in the quality of
this paper.

Furthermore, the author is grateful to his advisor, David J. Green,
for help, constructive comments and the interest in group cohomology.

\bibliography{library}

\end{document}